\documentclass[smallextended,nospthms]{svjour3}

\usepackage{graphicx}%
\usepackage{multirow}%
\usepackage{amsmath,amssymb,amsfonts}%
\usepackage{amsthm}%
\usepackage{mathrsfs}%
\usepackage[title]{appendix}%
\usepackage{xcolor}%
\usepackage{textcomp}%
\usepackage{manyfoot}%
\usepackage{booktabs}%
\usepackage{algorithm}%
\usepackage{algorithmicx}%
\usepackage{algpseudocode}%
\usepackage{listings}%

\usepackage{lineno}
\usepackage{comment}
\usepackage{subfigure}
\newtheorem{Theorem}{Theorem}[section]
\newtheorem{Example}{Example}[section]
\newtheorem{Lemma}{Lemma}[section]
\newtheorem{Definition}{Definition}[section]
\newcommand\mybox{\hbox to 0pt{}\hfill$\rlap{$\sqcap$}\sqcup$}

\begin{document}
	
	\title{Convex cones, assessment functions, balanced attributes}
	
	\author{Ignacy Kaliszewski}
	
	\authorrunning{Ignacy Kaliszewski} 
	
	\institute{Ignacy Kaliszewski \at Systems Research Institute, Polish Academy of Sciences,\\
		WIT University, Warsaw, Poland.\\
		\email{ignacy.kaliszewski@ibspan.waw.pl}     
			}

	\date{Received: date / Accepted: date}
	
	\maketitle
	
	\baselineskip = 1.5\baselineskip
	
	

\abstract{
 We investigate a class of polyhedral convex cones, with $R^k_+$ (the nonegative orthant in $\mathbb{R}^k$) as a special case. We start with the observation that for convex cones contained in $\mathbb{R}^k$, the respective cone efficiency is inconsistent with the Pareto efficiency, the latter being deeply rooted in economics, the decision theory, and the multiobjective optimization theory. Despite that, we argue that convex cones contained in $\mathbb{R}^k$ and the respective cone efficiency are also relevant to these domains.  
 To demonstrate this, we interpret polyhedral convex cones of the investigated class in terms of assessment functions, i.e., functions that aggregate multiple numerical attributes into single numbers.
 Further, we observe that all assessment functions in the current use share the same limitation; that is, they do not take explicitly into account attribute proportionality. In consequence, the issue of {\em attribute balance} (meaning {\it the balance of attribute values}) escapes them. In contrast, assessment functions defined by polyhedral convex cones of the investigated class, contained in $\mathbb{R}^k$, enforce the attribute balance. However, enforcing the attribute balance is, in general, inconsistent with the well-established paradigm of Pareto efficiency. We give a practical example where such inconsistency is meaningful.
} 
 
   
\section{Introduction}

We start with revisiting the notion of $K$-efficiency, where $K$ is a convex cone in $\mathbb{R}^k$.

\begin{Definition}
	\label{definition1.1}
	Given set $Z \subset \mathbb{R}^k$ and a closed convex cone $K$, element $\bar{y} \in Z$ is $K$-efficient, if there is no other element $y \in Z$ such that $y - \bar{y} \in K$.
\end{Definition}
 

The notion of $K$-efficiency, closely related to partial orders induced by convex cones, has been extensively studied in vector optimization since the fiftieth of the last century (cf. Kuhn, Tucker 1951, Yu 1974, Wierzbicki 1977, Henig 1982, Pascoletti and Serafini 1984, Jahn 1986, Luc 1989, Kaliszewski 1994). The interest in $K$-efficiency can be ascribed to the fact that it extends the notion of $R^k_+$-efficiency ({\it Pareto efficiency}), where $R^k_+ = \{ y  \in \mathbb{R}^k \, | \, y_l \geq 0, \ l = 1,\dots,k \}$, steaming from economics (Pareto 2014). To retain the relationship of $K$-efficiency to the economic inspiration and to keep its relevance for the {\it general decision theory}, the assumption $R^k_+ \subseteq K$ is adopted. With this assumption in force, if an element $\bar{y}$ of some set  $Z, \ Z \in \mathbb{R}^k$, is $K$-efficient, this element is also $R^k_+$-efficient. In other words, for any convex cone $K, \ R^k_+ \subseteq K$, $\bar{y}$ is not $K$-efficient if there exists $y \neq \bar{y}, \ \bar{y}_l \leq y_l, \ l = 1, \dots, k$. However, this statement is not in general valid for $K \subset R^k_+$. In consequence, when $K \subset R^k_+$, elements that are $K$-efficient are not necessarily $R^k_+$-efficient (Pareto efficient); in other words, they need not to adhere to the concept of Pareto efficiency.

Within the $R^k_+$-efficiency framework, many results are available for numerically deriving $R^k_+$-efficient elements. Here, one can mention:  

- Characterizations (i.e., necessary and sufficient conditions) of (properly) $R^k_+$-efficient elements under the assumption that set $Z$ is convex (Klinger 1967, Geoffrion 1968); such characterizations are based on optimizing the weighted sum of element components.

- Characterizations of (properly) $R^k_+$ efficient elements under no assumptions about set $Z$ except that it is bounded
(Choo, Atkins 1983, Steuer, Choo 1983, Steuer 1986, Kaliszewski 1987,1994, Wierzbicki 1999); such characterizations are based on the Leontief (in the multiple criteria decision-making and the multiobjective optimization -- the Chebyshev) function defined over elements $y \in Z$.   

Such characterizations apply only to a subset of $R^k_+$-efficient elements of a set, namely to properly $R^k_+$-efficient elements. This fact was observed first by Kuhn and Tucker (1951) in the context of optimizing the weighted sums of components of elements $y \in Z$. A formal definition of properly $R^k_+$-efficient elements was given by Geoffrion (1968). The Geoffrion definition of proper efficiency was subsequently generalized to arbitrary convex cones (cf. Wierzbicki 1977, Borwein 1977, Benson 1979, Nieuwenhuis 1981, Henig 1982). However, the definition by Geoffrion is the only quantitative one, whereas the others are qualitative. In this work, we confine ourselves to the definition by Geoffrion. Nonetheless, in $R^k_+, \ k < +\infty$ (i.e., in finite-dimensional real spaces) all those definitions are equivalent (Henig 1982).

As shown below, characterizations based on the Leontief function apply also to the investigated class of polyhedral convex cones $K$ when $R^k_+ \subset K$. Thus, they also provide for the effective numerical derivation of $K$-efficient elements for such cones. We also show that such characterizations work for the case $K \subset R^k_+$.
For the case $K \subset R^k_+$, to our best knowledge no effective characterizations of $K$-efficient elements exist, though the possibility of deriving such elements with extensions of the Leontieff function was considered in Kaliszewski (2006). The reason for the lack of interest in that sort of results may be explained by that researchers perceive no need to depart from the Pareto equilibrium paradigm, which is so deeply rooted in the economic and decision-making theories (cf., e.g., Fudenberg, Tirole 1991, Chinchuluun et al. 2008, Bruni 2013). 

When working with cones $K$ such that $K \subset R^k_+$, one have to accept that for a pair of $K$-efficient elements $y, y' \in Z, \ y \neq y'$, the relation $y' \neq y, \ y_l \leq y'_l, \ l = 1, \dots, k$, can hold. In other words, elements that are $K$-efficient, $K \subset R^k_+$,  are not necessarily $R^k_+$-efficient. In the context of the general decision theory, this is counter-intuitive. In this paper, we point to  specific circumstances where $K$-efficient elements, $K \subset R^k_+$, that are not $R^k_+$-efficient, are meaningful.

The outline of the paper is as follows. In Section~\ref{intro}, we recall basic definitions and the necessary notation. There we also define or rather recall, a~class of polyhedral cones, with $R^k_+$ as the special case. In Section~\ref{characterization}, we recall Geoffrion's definition of properly $R^k_+$-efficiency and the characterization of such elements. A generalization of Geoffrion's definition for polyhedral convex cones is given in Section~\ref{poly}. This generalization applies to convex cones $K$ that contain $R^k_+$ (the case where $K$-efficient elements are $R^k_+$-efficient), as well as to those that are contained in $R^k_+$ (the case where $K$-efficient elements not necessarily are $R^k_+$-efficient). In Section~\ref{nP_rankings}, the relevance of the presented results to domains that make use of quantitative multiattribute evaluations, and to rankings in particular, is discussed. The effect of ranking with different aggregating functions is illustrated in the example of the global Environmental Performance Index. Section~\ref{conclusions} concludes.

\section{Basic definitions and notation}
\label{intro}

\begin{figure}
	\centering
	\includegraphics[scale=0.4]{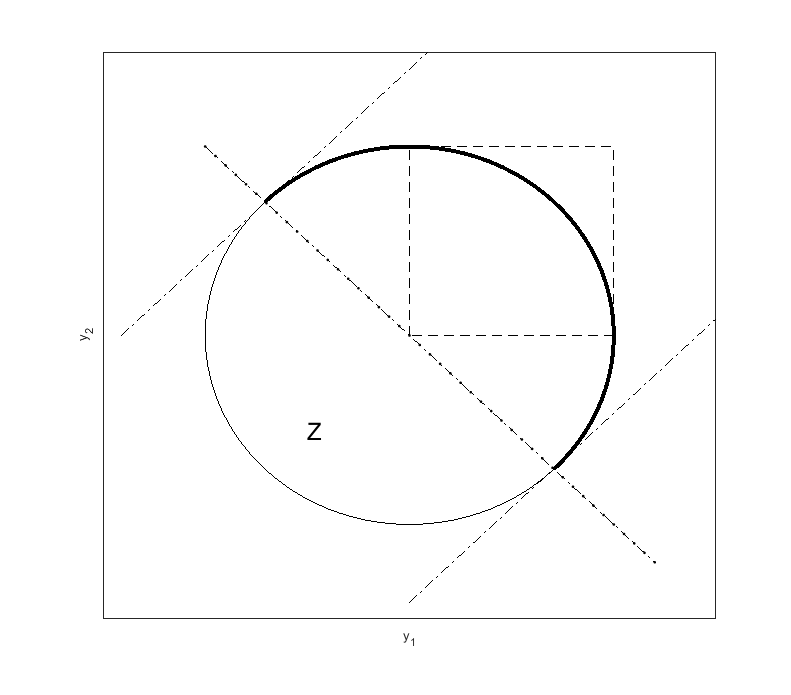}
	\caption{The thick curve indicates elements that can be $K$-efficient under an appropriate selection of a convex cone $K$. In particular, all of them are $K^{Poly}_{-\frac{1}{k}}$-efficient. From those elements, all elements that belong to the dotted square are $R^k_+$-efficient, whereas the remaining ones are not $R^k_+$-efficient (are not Pareto efficient). The latter are $K^{Poly}_{\rho}$-efficient for some $-\frac{1}{k} \leq \rho < 0$.} 
	\label{Figure1}
\end{figure}

	In $\mathbb{R}^{k}$ we consider polyhedral cones 
	\begin{equation}
		\label{Kpoly}
	K^{Poly}_\rho = \{ y \ | \ y_l + \rho \sum^k_{i=1} y_i \geq 0, \ l = 1,\dots, k\}, \ -\frac{1}{k} \leq \rho < +\infty,
	\end{equation} 
	where $\rho = +\infty$ means $\rho \rightarrow +\infty$. $K^{Poly}_\rho$ cones are closed and convex. In $\mathbb{R}^2$, convex cones are polyhedral.
	
	For $K^{Poly}_\rho$ cones, Definition \ref{definition1.1} can be restated as follows.

\begin{Definition}
	\label{definition2.1}
	Given set $Z$ and $K^{Poly}_\rho$ cone, element $\bar{y} \in Z$ is $K^{Poly}_\rho$-efficient, if there is no element $y$ in $Z$, $y \neq \bar{y}$, such that $y_l + \rho \sum^k_{i=1} y_i \geq \bar{y}_l + \rho \sum^k_{i=1} \bar{y}_i, \ l = 1,\dots,k$.
\end{Definition}

To see the equivalence of this two definitions when $K = K^{Poly}_\rho$, suppose first  that $y - \bar{y} \in K^{Poly}_\rho$. Then, $(y_l - \bar{y}_l) + \rho \sum^k_{i=1} (y_i - \bar{y}_i) \geq 0$, $l = 1,\dots,k$, and by the Definition \ref{definition2.1} $\bar{y}$ is not $K^{Poly}_\rho$-efficient.

Suppose now that $y_l + \rho \sum^k_{i=1} y_i \geq \bar{y}_l + \rho \sum^k_{i=1} \bar{y}_i, \ l = 1,\dots,k$. Then, $(y_l - \bar{y}_l) + \rho \sum^k_{i=1} (y_i - \bar{y}_i) \geq 0$,  $l = 1,\dots,k$, thus $y - \bar{y} \in K^{Poly}_\rho$ and by the Definition \ref{definition1.1}, $\bar{y}$ is not $K^{Poly}_\rho$-efficient.

\vspace{0.2cm}
For $\rho = 0$, $K^{Poly}_\rho$ cones reduce to $R^k_+ = \{ y \ | \ y_l \geq 0, \ l = 1,\dots, k\}$. 

For $\rho = +\infty$, $K^{Poly}_\rho$ cones reduce to hyperplane $H_+ = \{ y \ | \ \sum^k_{i=1}y_i \geq 0 \}$. 

For $\rho = -\frac{1}{k}$, $K^{Poly}_\rho$ cones reduce to half line $\vec{1}t, \ t \in [0:+\infty)$, where $\vec{1}$ denotes $k$-dimensional vector of ones. To see this, observe that all inequalities $y_l -\frac{1}{k}\sum^k_{i=1} y_i \geq 0, \ l = 1,\dots,k$, where the second term of each inequality is the mean value of all components, hold only if all $y_l$ are equal, i.e., they are located on $\vec{1}t, \ t \in [0:+\infty)$. 

Figure \ref{Figure1} shows all elements that can be efficient under an appropriate selection of $K^{Poly}_\rho$ cones. All elements of the thick curve are $K^{Poly}_{-\frac{1}{k}}$-efficient. Set $Z$, represented in this figure schematically by a ball, can be of any nature, e.g., it can be disconnected (like two disjoint balls) or discrete.

\begin{figure}
	\centering
	\includegraphics[scale=0.4]{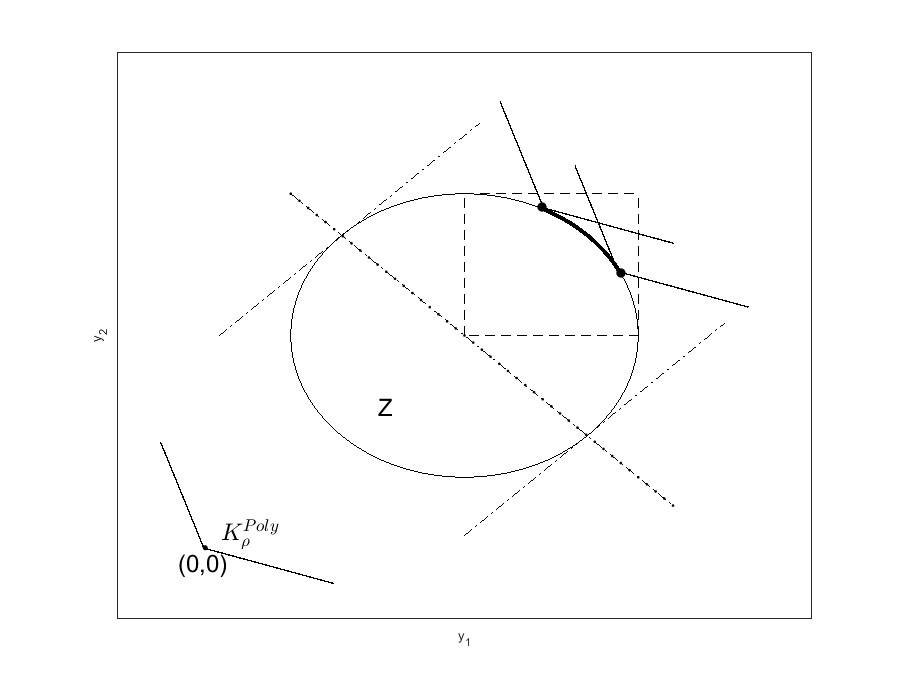}
	\caption{The thick curve indicates elements that  are $K^{Poly}_{\rho}$-efficient for some $\rho > 0$.} 
	\label{Figure2}
\end{figure}

Figure \ref{Figure2} presents $K^{Poly}_{\rho}$-efficient elements for some $\rho > 0$.

\section{Properly $R^k_+$-efficient elements}
\label{characterization}

In this section, we recall the definition of properly $R^k_+$-efficient elements by Geoffrion and a numerically tractable characterization of such elements. The characterization provides for the numerical derivation of properly $R^k_+$-efficient elements. We provide a generalization of this definition to the class of polyhedral convex cones under investigation. This generalization applies to $K^{Poly}_{\rho}$ cones that contain $R^k_+$, as well as to $K^{Poly}_{\rho}$ cones that are contained in $R^k_+$. We provide also a numerically tractable characterization of properly $K^{Poly}_{\rho}$-efficient elements. 

In the following definition the terms {\em efficient} and {\em properly efficient} refer to  $R^k_+$-efficiency and proper $R^k_+$-efficiency.

\begin{Definition}
	\label{definition3.1}{(Geoffrion 1968)}
	An element $\bar{y} \in Z$ is properly  efficient if it is efficient and there exists a finite number $M > 0$ such that for each $l, \ l=1,\dots,k$,
	we have
	\[
	\frac{y_l - \bar{y}_l}{\bar{y}_j - y_j} \leq M
	\]
	for some $j$ such that
	$y_j < \bar{y}_j \, ,$ whenever $y \in Z$ and
	$y_l > \bar{y}_l \, .$
\end{Definition}

The following theorem gives a characterization of properly $R^k_+$-efficient elements of $Z$. 

By $R^k_>$ we denote $\{ y \ | \ y_l > 0, \ l = 1,\dots, k\}$.

From now on we assume that $Z \in R^k_>$ (all components of any $y \in Z$ are positive, cf. Section \ref{poly}). With that assumption in force, we have the following theorem.

\begin{Theorem}
	\label{theorem3.1}
An element
	$\bar{y}\in Z$ is properly $R^{k}_+$-efficient if  and  only  if
	there exists 
	$\lambda \in R^k_> $ and 
	$\rho>0 \, , $ such that  $\bar{y}$  solves
	\[
	\begin{array}{ll}
		P^{\infty}_{R^{k}_+} \, : \index{problem!$P^\infty$} \ \ \ \ \ \ \ \ \ \ \ \ &  \max_{y\in Z}\min_{1\leq l \leq k}
		\lambda_{l}(y_{l} + \rho \sum^k_{i=1} y_i) \, .
	\end{array}
	\]
	
	For each properly efficient element
	$y\in Z$  there exists
	$\lambda \in R^k_>$
	such  that  $y$  solves
	$P^{\infty}_{R^{k}_+}$   uniquely for every
	$\rho>0$  satisfying
	$M\leq((k-1)\rho)^{-1} \, .$
	
	For each element  $y \in Z$ which solves
	$P^{\infty}_{R^{k}_+}$, the inequality
	$M\leq(1+(k-1)\rho)\rho^{-1}$ holds.
\end{Theorem}

The proof of this theorem follows from the proof of Theorem \ref{theorem4.1} as a special case.

This theorem is similar to one given in Kaliszewski 1994 (Theorem 4.2, pp 48-50), which in turn was based on the original formulation and proof given in Choo and Atkins, 1983. 

\section{Properly $K^{Poly}_\rho$-efficient elements}
\label{poly}

Definition \ref{definition3.1} can be generalized to  $K^{Poly}_\rho$ cones.

Denote $r_l = y_l + \rho \sum^k_{i=1} y_i, \ y \in Z, \ l = 1,\dots,k, \ r = (r_1,\dots,r_k)$, where $-\frac{1}{k} \leq \rho < +\infty$. 

The following definition is just the Definition \ref{definition3.1} with components $y_l$ defining  $R^k_+$ replaced by components $r_l, \ l=1,\dots,k$, defining $K^{Poly}_{\rho}$ cones (\ref{Kpoly}). 

\begin{Definition}
	\label{definition4.1}
	An element $\bar{y} \in Z$  is  properly   $K^{Poly}_{\rho}$- efficient if it is $K^{Poly}_{\rho}$-efficient and there exists a finite number $N > 0$ such that for each $l, \ l=1,\dots,k$,
	we have
	\[
	\frac{r_l - \bar{r}_l}{\bar{r}_j - r_j} \leq N
	\]
	for some $j$ such that
	$r_j < \bar{r}_j \, ,$ whenever $r \in Z$ and
	$r_l > \bar{r}_l \, .$
\end{Definition}

\begin{Example}
	\label{example}
	
	Let $Z = \{ y \ | \ y^2_1 + y^2_2 \leq 1, \ y_1, y_2 \geq 0 \}$, $\rho = 1$. Then, $K^{Poly}_{\rho} = \{y \ | \ r_1 = y_1 + \rho(y_1 + y_2) \geq 0, \ r_2 = y_2 + \rho(y_1 + y_2) \geq 0 \}$. 
	
	The $R^k_+$-efficient element $\bar{y} \in Z$, where the line $ y_1 + \rho(y_1 + y_2) = \sqrt{5}$ is tangent to $Z$, has coordinates $\bar{y}_1 = \frac{2\sqrt{5}}{5}$, $\bar{y}_2 = \frac{\sqrt{5}}{5}$. Indeed, $(\frac{2\sqrt{5}}{5})^2 + (\frac{\sqrt{5}}{5})^2 = 1$. Element $\bar{y}$ is not properly (is {\upshape improperly}) $K^{Poly}_\rho$-efficient. This is because for a sequence of elements on the circle $y^2_1 + y^2_2 = 1$, tending to element $\bar{y}$ from the right side (i.e., with decreasing values of $y_1$), the sequence of corresponding numbers $\frac{r_2 - \bar{r}_2}{\bar{r}_1 - r_1}$ is unbounded. The proof of this fact is given in Appendix 1.
\end{Example}

Figure \ref{Figure3} presents examples of elements that are improperly  $K^{Poly}_{\rho}$-efficient for $\rho =0$ (i.e. for $K^{Poly}_\rho = R^k_+$), and for some $\rho < 0$.

\begin{figure}
	\centering
	\includegraphics[scale=0.4]{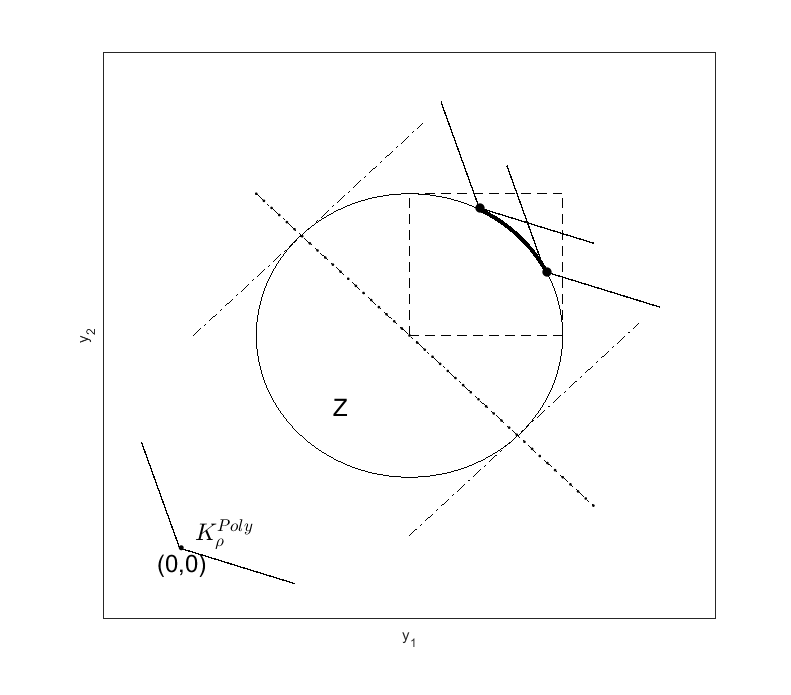}
	\caption{Elements: improperly $R^k_+$-efficient (circles) and improperly $K^{Poly}_{\rho}$-efficient (discs) for some $\rho < 0$.} 
	\label{Figure3}
\end{figure}

With Definition \ref{definition4.1}, we can generalize Theorem \ref{theorem3.1} to Theorem \ref{theorem4.1}.

The crucial fact for the correctness of Theorem \ref{theorem3.1} is that all $r_l = y_l + \rho \sum^k_{i=1} y_i, \ y \in Z, \ l = 1,\dots,k$, are positive. This is the consequence of the assumptions $Z \subset 
R^k_>$ and $\rho \geq 0$. However, admitting $\rho < 0$, this is no longer guaranteed. But as $Z$ is bounded by assumption, we can offset (translate) set $Z$ by $r^\bullet = (r^\bullet_1,\dots,r^\bullet_k)$, where $r^\bullet_l = min_{r_l \in Z} \, r_l + \varepsilon$, $\varepsilon > 0$. Thus, $Z^\bullet = \{ r^\bullet \} + Z$ and by this for any $y \in Z^\bullet$ the respective $r_l, \ l=1,\dots,k$, are positive $(Z^\bullet \subset R^k_>$).

\begin{Theorem}
	\label{theorem4.1} 
	Given $\rho$, $-\frac{1}{k} \leq \rho < +\infty$\, an element
	$\bar{y} \in Z^\bullet$ is properly $K^{Poly}_{\rho}$-efficient if and only if
	there exists
	$\lambda \in R^k_>$ and 
	$\sigma > 0$, such that  $\bar{y}$  solves
	\[
	\begin{array}{ll}
		P^{\infty}_{K^{Poly}_{\rho}} \, :  \ \ \ \ \ \ \ \ \ \ &  \max_{y\in Z^\bullet} \min_{1\leq l \leq k}
		\lambda_{l}(r_{l} + \sigma \sum^k_{i=1} r_i) \, .
	\end{array}
	\]
	
	For each properly $K^{Poly}_\rho$-efficient element
	$y\in Z^\bullet$  there exists 
	$\lambda \in R^k_>$
	such  that  $y$  solves
	$P^{\infty}_{K^{Poly}_{\rho}}$   uniquely for every
	$\sigma > 0$  satisfying
	$N\leq((k-1)\sigma)^{-1} \, .$
	
	For each element  $y \in Z^\bullet$   which solves
	$P^{\infty}_{K^{Poly}_{\rho}}$  the inequality
	$N\leq(1+(k-1)\sigma)\sigma^{-1}$ holds.
\end{Theorem}

The proof is given in Appendix 2.

The following lemma is the immediate consequence of Theorem \ref{theorem4.1}.

\begin{Lemma}
	\label{theorem4.2} 
	Given $-\frac{1}{k} \leq \rho < +\infty$, an element
	$\bar{y}\in Z^\bullet$ is $K^{Poly}_{\rho}$-efficient only if there exists  
	$\lambda \in R^k_>$ such that  $\bar{y}$ uniquely solves
	\[
	\begin{array}{ll}
		\bar{P}^{\infty}_{K^{Poly}_{\rho}} \, :  \ \ \ \ \ \ \ \ \ \ &  \max_{r\in Z^\bullet}\min_{1\leq l \leq k}
		\lambda_{l} r_{l} \, .
	\end{array}
	\]
\end{Lemma}

The proof of this lemma follows directly from the necessity part of the proof of Theorem \ref{theorem4.1} with parameter $\sigma$ set to $0$.  

The proof of the necessity part of Lemma \ref{theorem4.2} is constructive, i.e., it shows how to construct $\bar{\lambda} \in R^k_>$ for which $P^{\infty}_{K^{Poly}_{\rho}}$ yields a given $K^{Poly}_{\rho}$-efficient element $\bar{y}$. This immediately provides a test for  $K^{Poly}_{\rho}$-efficiency, in particular for $R^k_+$-efficiency, (cf. Kaliszewski et al. 2016). The test goes as follows.

\vspace{0.2cm}
\noindent
{\bf ${\mathbf K^{Poly}_{\rho}}$-efficiency test}

\vspace{-0.3cm}
\begin{enumerate}
	\item For element $\bar{y} \in Z^\bullet$ tested for $K^{Poly}_{\rho}$-efficiency, set
	\[
	\bar{\lambda}_l =  (\bar{r}_{l})^{-1}, \ \ l = 1,\dots,k,
	\]
	where $\bar{r}_l = \bar{y}_l + \rho \sum^k_{i=1} \bar{y}_i$.
	\item Solve problem $P^{\infty}_{K^{Poly}_{\rho}}$ with $\bar{\lambda}$ and $\sigma = 0$.
	
	\item If $\bar{y}$ solves $P^{\infty}_{K^{Poly}_{\rho}}$ with  $\bar{\lambda}$ and $\sigma = 0$, it is  $K^{Poly}_{\rho}$-efficient, otherwise it is not $K^{Poly}_{\rho}$-efficient.
\end{enumerate}

\begin{figure}
	\centering
	\includegraphics[scale=0.4]{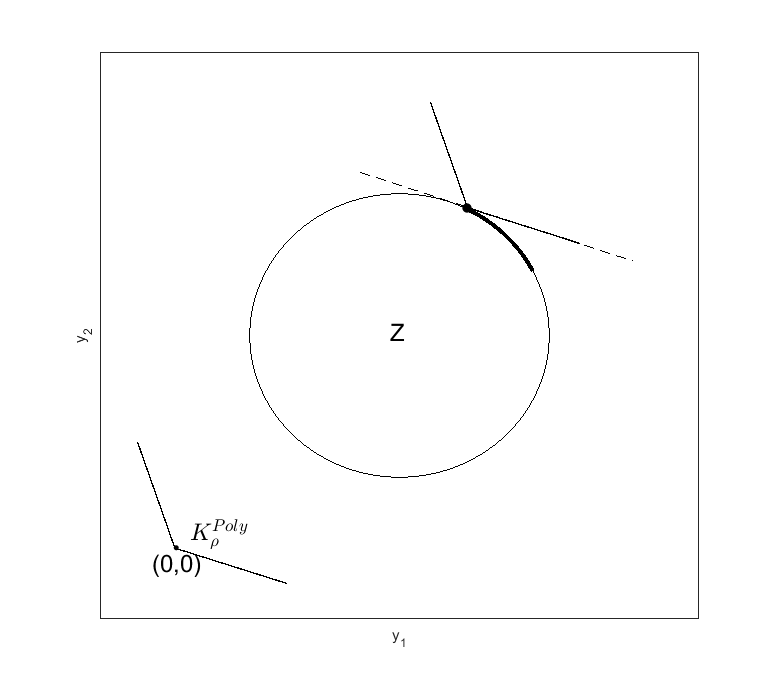}
	\caption{An illustration to Example \ref{example}: element $\bullet$ is improperly $K^{Poly}_{\rho}$-efficient.} 
	\label{Figure4}
\end{figure}

\section{Assessment functions}
\label{assess_functions}

Each domain of human activities faces the challenge of deciding which action or object, collectively named {\it alternatives}, in a pool is the best, and often also what is the next best, the third best, and so on. Whenever possible, time is not critical, and skills are available, analytic methods are applied to choose the best alternative according to a predefined meaning of being {\it best}.

Alternatives are characterized by their attribute values. To facilitate analyses, it is a common practice to map qualitative attribute values into numbers. If alternatives are represented by only one attribute, the choice of the best alternative reduces to sorting and selecting the alternative (and the following best, and so on) with the attribute’s best value. Otherwise, multiple attributes induce on the set of alternatives a partial order, and the best alternative can no longer be chosen by sorting. Various approaches, generally termed {\em rankings}, have been proposed to aggregate multiple attributes into one-dimensional surrogate measures. However, this amounts to replacing partial orders, spanned by multiple attributes, with linear orders. Except in the case they coincide, this cannot be done but with a dose of arbitrariness.

If alternatives' attributes $y$ are not given explicitly but by functional mappings (like in, e.g., multiobjective optimization), considerations are presented in terms of both $x$ and $y$, where $y = f(x), \ Z = f(X_0), \ f = (f_1,\dots,f_k), \ f_l: R^n \rightarrow R, \ l=1,\dots,k$, and $x$ belong to some set of alternatives $X_0 \subset \mathbb{R}^n$. 


The standard way to reduce partial orders to linear orders is via multivariate {\em assessment functions}. Depending on the adopted assessment function, different rankings result. Varying by application domains, multiattribute assessment functions are termed ''value functions'', ''utility functions'',  ''production functions'', or ''scalarizing functions''. 

For consistency and generality of the presentation, below we assume that alternative attribute values are given by elements of set $Z \subset R^k$.
The standard assumption in application-oriented domains, like the production theory or the utility theory, is that attribute values are positive, i.e., $Z \subseteq R^k_>$, where $R^k_> = \{ y \ | \ y_l > 0, \ l = 1,\dots,k \}$ and we follow this convention. 

Consider the general form of assessment functions:
\begin{equation}
	\label{eq_asses_funct_1}
	{\cal F}(y),
\end{equation}
where $y = (y_1, \dots, y_k)$, \ $y_l, \ l=1,\dots,k$, are attribute values, ${\cal F}$ is a strictly increasing function, i.e., $y_l \leq y_l$, $l=1,\dots,k$, and $y_l < y_l'$ for some $l$, implies ${\cal F}(y) < {\cal F}(y')$, or at least it is an increasing function, i.e., $y_l \leq y_l'$, $l=1,\dots,k$, implies ${\cal F}(y) \leq {\cal F}(y')$.

To accommodate the observation that people prefer increases in value (wealth) less and less as the value increases, assessment functions 
are assumed to be concave (cf., e.g.,   Bourgeois-Gironde 2020).

For presentation clarity, below we argue as if all attributes are of equal importance. To take into account situations when some attributes are relatively less or more important than others, assessment functions admit positive {\em weights} for attribute weighting. Weighted attributes can be considered as new attributes with no weights, precisely, with all weights equal to $1$.

\subsection{Assessment functions in the production theory and the utility theory}
\label{prod_theory}

The classical assessment production functions and utility functions are based on the mean of the order $p$: 
\begin{equation}
	\label{eq_asses_funct_4}
	{\cal F}(y) = (\sum_{l=1}^k y_l^p)^\frac{1}{p}. 		 		
\end{equation}

For $p \rightarrow 0, \ {\cal F}(y) = \prod_{l=1}^k y_l$, and for $p\rightarrow -\infty$, ${\cal F}(y) = \min_{l, \, 1 \leq l \leq k} y_l$ (the Leontief function).

Except for the Leontief function, functions (\ref{eq_asses_funct_4}) are strictly increasing and for $p \leq 1$ they are concave. The Leontief function is not strictly increasing but only increasing.

In the production theory, the parameterized version of functions (\ref{eq_asses_funct_4}) for $-\infty < p < 1$ are known as the (family of) CES functions (cf., e.g., de La Grandville 2016):

\[
{\cal F}(y) = a(\sum^k_{l=1} a_l y_l^p)^{\frac{1}{p}}, \ p < 1,
\]
and for $p \rightarrow 0$, as the Cobb-Douglass production function: 

\[ 
{\cal F}(y) = a\prod^k_{l=1} y_l^{\alpha_l}, 
\]
(cf., e.g., Jacques 2018). For $p =1$, (\ref{eq_asses_funct_4}) is the linear function.

\subsection{Pareto dominance-consistent assessment functions}
\label{PDCA}

If an alternative is not $R^k_+$-efficient in a given set of alternatives, (i.e., it is not Pareto efficient), then in this set there exists another alternative with all attribute values greater or equal and at least one attribute value strictly greater. We say that the latter alternative $R^k_+$-{\it dominates} the former, and the former is $R^k_+$-{\it dominated} by the latter. 

All functions (\ref{eq_asses_funct_4}) attain their maximal values over a set of alternative attribute values at the $R^k_+$-efficient alternatives\footnote{\mbox{ }In the case of the Leontief assessment function, this is true only if such an alternative is unique. If it is not unique, at least one of such alternatives is Pareto efficient. This is because the Leontieff assessment function is increasing but not strictly increasing (cf. Kaliszewski 1994, Miettinen 1999, Ehrgott 2005).} and  $R^k_+$-dominated alternatives yield lower value of the assessment functions than $R^k_+$-dominating alternatives. 
We call that property the {\it Pareto dominance consistency}, and functions with that property the {\it Pareto dominance-consistent assessment functions} (PDCA functions). 

Of special interest are two assessment functions. 

The Leontief assessment function, denoted ${\cal F}^L(y)$ (in the production theory, it is known as the Leontief production function, cf., e.g., Miller and Blair, 2009):

\begin{equation}
	\label{Leontief_function}
	{\cal F}^L(y) = \min_{l, \, 1 \leq l \leq k}y_l.
\end{equation} 
It takes the zero value at $\vec{0} = (0,\dots,0)$. As the limit case od (\ref{eq_asses_funct_4}), it is also a~PDCA function.

Another example of a PDCA function is the Chebyshev assessment denoted ${\cal F}^{Ch}(y)$ function, widely used in multiobjective optimization and multiple criteria decision-making (cf., e.g., Wierzbicki 1977,1999, Ehrgott et al. 2004, Ehrgot 2005, Miettinen 1999, Kaliszewski et al. 2016), namely
\begin{equation}
	\label{Chebyshev}
	{\cal F}^{Ch}(y) = \min_{l, \, 1 \leq l \leq k}(y_l - y_l^{ref}) = - \max_{l, \, 1 \leq l \leq k}(y_l^{ref} - y_l),
\end{equation} where $y^{ref}$, called {\it reference point}, is any element of $\mathbb{R}^k$. If $y_1^{ref} = \dots = y_k^{ref} = a$, and $a > \max_{l, \, 1 \leq l \leq k}(y_l^{ref} - y_l)$, then ${\cal F}^L(y) = {\cal F}^{Ch}(y) - a$. For $y^{ref} = \vec{0}$, the Leontief assessment function and the Chebyshev assessment function coincide. 

\begin{figure}[h]
	\centering
		\subfigure[]{\includegraphics[width=0.45\textwidth]{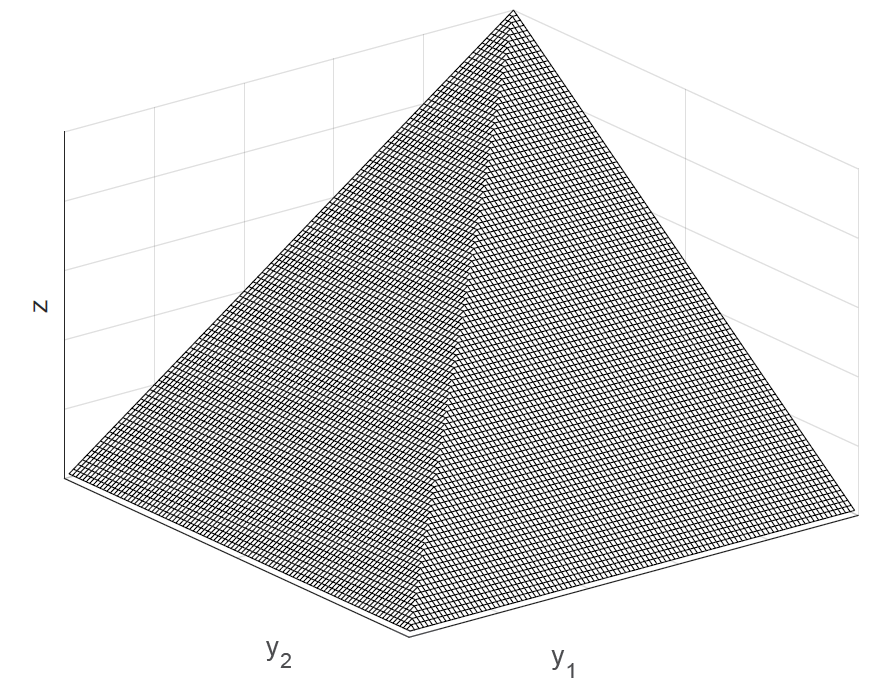}}
		\subfigure[]{\includegraphics[width=0.45\textwidth]{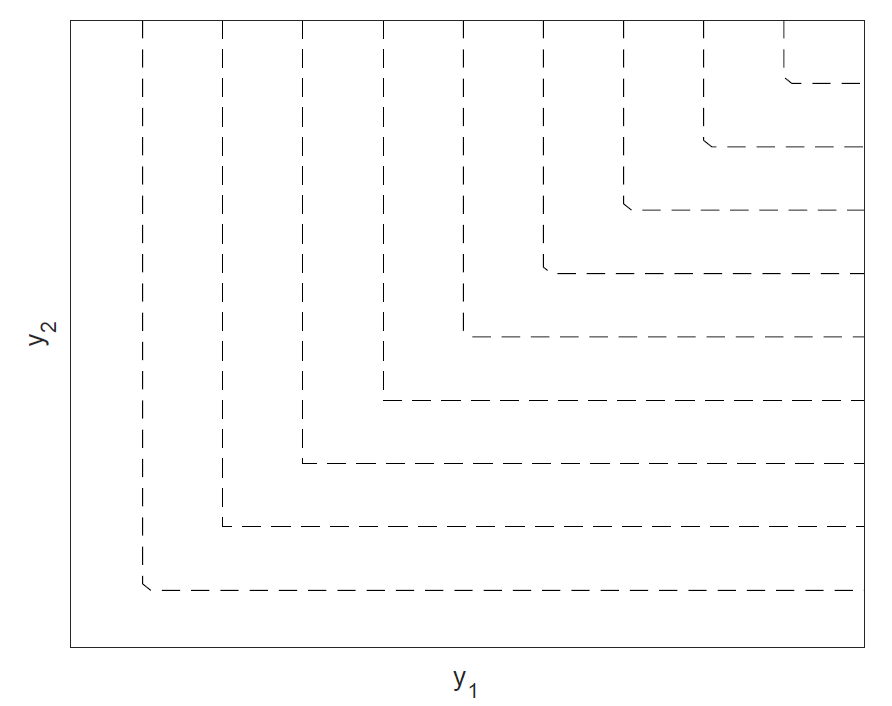}}
	\caption{Graph and contours of the Leontief function (\ref{Leontief_function}) in $\mathbb{R}^2$; a PDCA function.}
	\label{Figure2mat}
\end{figure}

Figure \ref{Figure2mat} presents the graph and contours of the Leontief function. The contours of the Leontief function (and also of the Chebyshev function) are shifted cones $R^k_+$ and that explains its particular relevance to the numerical derivation of $R^k_+$-efficient elements. These two functions can be interpreted as  analytical tools to verify whether an element is $R^k_+$-efficient by applying directly Definition \ref{definition2.1}.

Yet another example of a PCDA function is the augmented Leontief function is defined as
\begin{equation}
	\label{augLeontief}
{\cal F}(y) = \min_{1 \leq l \leq k} (y_l + \ \rho\sum_{i = 1}^{k}y_i), \ 0 \leq \rho < +\infty.
\end{equation}
For $ \rho > 0$ this function constitutes the objective function of optimization problem $P^{\infty}_{R^{k}_+}$ in Theorem \ref{theorem3.1}. At the limit $\rho \rightarrow + \infty$, function (\ref{augLeontief}) reduces to the linear function, which also is a PCDA function.

\begin{figure}
	\centering
		\subfigure[]{\includegraphics[width=0.45\textwidth]{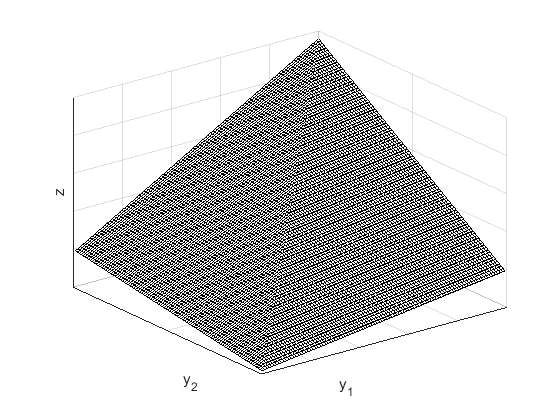}}
		\subfigure[]{\includegraphics[width=0.45\textwidth]{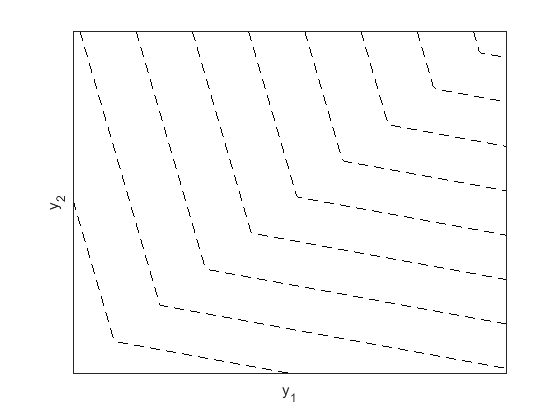}}
	\caption{Graph and contours of the Generalized Leontief function (\ref{GLF}) in $\mathbb{R}^2$ for some $\rho, \ 0 < \rho < + \infty  $; a PDCA function.}
	\label{Figure3mat}
\end{figure}

Figure \ref{Figure3mat} presents the graph and contours of the augmented Leontief function. The contours are shifted cones $K^{Poly}_{\rho}$ for some $\rho > 0$.

PDCA functions are increasing functions, so if some attribute values increase and the remaining are kept constant, the functions assume greater or equal values, irrespective of the attribute balance (proportions). When the required attribute balance should be maintained, or at least the deviation from it should be minimized (cf. Section \ref{nP_rankings} for a practical interpretation of such a requirement), assessment functions that take into account these effects are needed. 

\subsection{Pareto dominance-inconsistent assessment functions}
\label{PDIA}

Alternatives with equal attribute values are {\it perfectly balanced}. Attributes of perfectly balanced alternatives are located on the half line $t \vec{1}, \ t \in [0:+\infty)$. Perfectly balanced alternatives represent the best alternatives in the sense of the attribute balance, the ''goodness'' of other alternatives depends on how much their attribute values deviate from the half line $t\vec{1}, \ t \in [0:+\infty)$. From two different elements of that half line, one $R^k_+$-dominates the other but since they are of equal ''goodness'' in the sense of the attribute balance, no PDCA function applies to that case. 

To assess alternatives concerning the degree they deviate from the attribute perfect balance we can use the  Generalized Leontief functions (with the augmented Leontef functions as the special case)
\begin{equation}
\label{GLF}
{\cal F}(y) = \min_{1 \leq l \leq k} (y_l + \ \rho\sum_{i = 1}^{k}y_i),	\  -\frac{1}{k} \leq \rho < +\infty,
\end{equation}
and select $\rho = -\frac{1}{k}$. As proved in Appendix 3, for $\rho = -\frac{1}{k}$ contours of the Generalized Leontief function reduce to parallel lines, as shown in Figure \ref{Figure6mat} in $\mathbb{R}^2$. The largest value of the function, namely zero, is represented by the contour defined by alternatives with all attribute values equal. This also holds for in $\mathbb{R}^k, \ k > 2$, as illustrated by Example \ref{example5.1} and Figure \ref{Figure5A} that shows points in $\mathbb{R}^3$ equidistant, in terms of function (\ref{GLF}) with $-\frac{1}{k}$, to the half line $ \vec{1}t, \ t \in [0:+\infty)$.

\begin{Example}
	\label{example5.1}
	Figure \ref{Figure9} shows attribute values located on a curve in $\mathbb{R}^3$. All the corresponding alternatives, except the one marked by the bullet, are $R^3_+$ dominated. For $\rho = \frac{1}{3}$, for the alternative $R^k_+$-dominated by all other alternatives represented on the curve (i.e., that one with attribute values on the half line $y = \vec{1}t, \ t \in [0:+\infty)$), function (\ref{augLeontief}) assumes value $0$, and for the other two elements it assumes 
	value $-3.3269$ (the element marked by circle) and value $-6.1568$ (the element marked by bullet).
	The curve is defined as
	
	\begin{equation}
		\begin{array}{c}
			t \in [-\frac{\pi}{4};0], \\ \\
			x = 6cos(t) - 2.2426, \\ \\
			y = 10sin(t) + 9.0710, \\ \\
			z = [2;7] \,.
		\end{array}
	\end{equation} 
\end{Example}
	
Depending on $\rho$, $ -\frac{1}{k} < \rho < +\infty$, the Generalized Leontief functions reflect the attribute balance to a different degree, the smaller $\rho$ the more the attribute balance comes to importance. However, only for $ 0 \leq \rho \leq +\infty$ these functions are PDCA functions, and for $ -\frac{1}{k} \leq \rho < 0 $ they are {\it Pareto dominance-inconsistent assessment functions} (PDIA functions).

\begin{figure}
	\centering
		\subfigure[]{\includegraphics[width=0.45\textwidth]{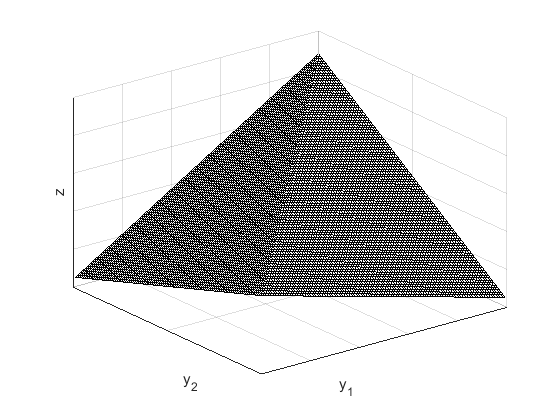}}
		\subfigure[]{\includegraphics[width=0.45\textwidth]{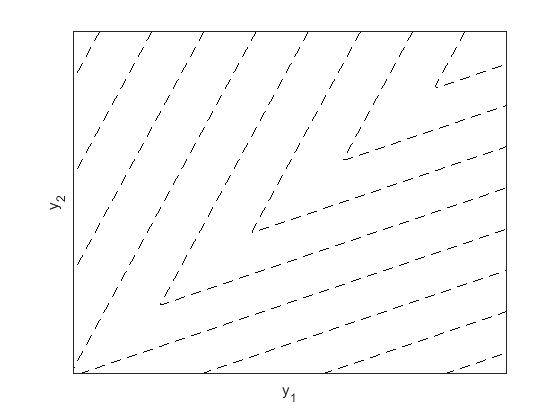}}
	\caption{Graph and contours of the Generalized Leontief function (\ref{GLF}) in $\mathbb{R}^2$ for some $\rho, \ -\frac{1}{2} < \rho < 0 $; a PDIA function.}
	\label{Figure5mat}
\end{figure}

\begin{figure}
	\centering
	\subfigure[]{\includegraphics[width=0.45\textwidth]{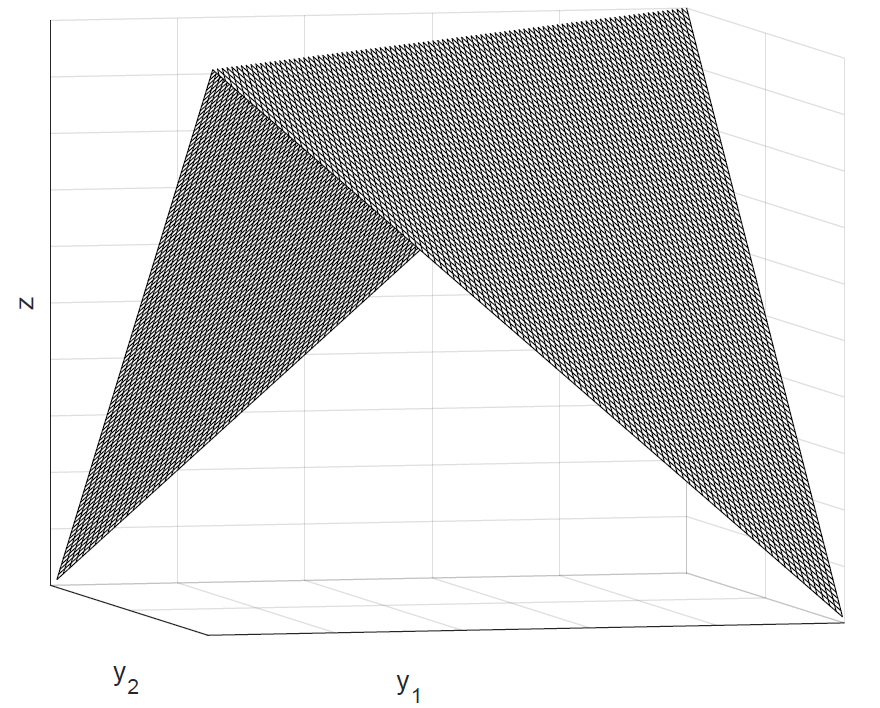}}
		\subfigure[]{\includegraphics[width=0.45\textwidth]{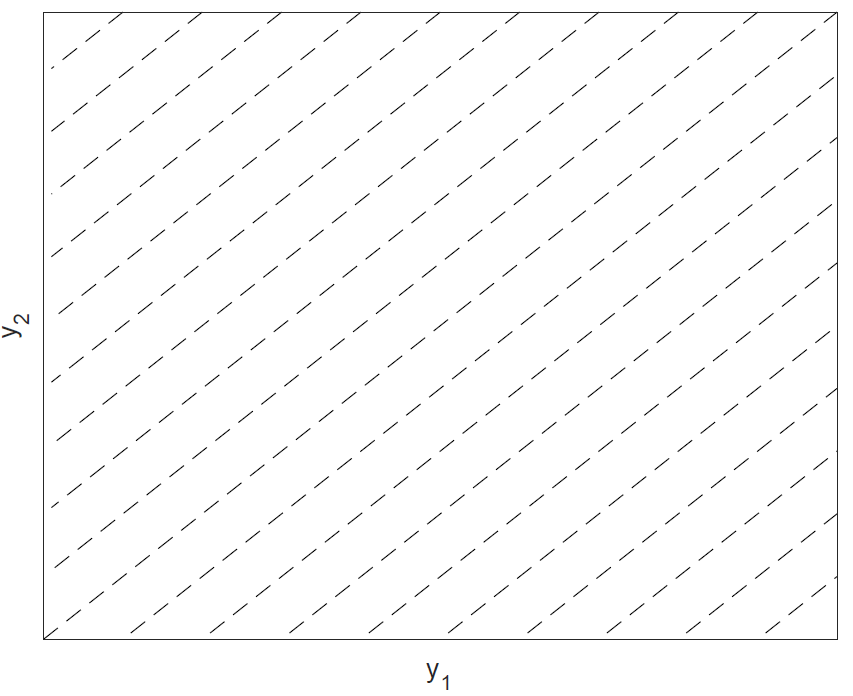}}
	\caption{Graph and contours of the Generalized Leontief function (\ref{GLF}) in $\mathbb{R}^2$  for $ \rho = -\frac{1}{2} $; a~PDIA function.}
	\label{Figure6mat}
\end{figure}

\begin{figure}
	\centering
	\includegraphics[scale=0.45]{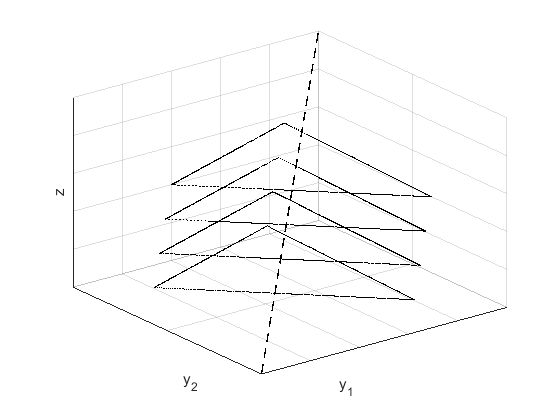}
	\caption{The equilateral triangles represent in $\mathbb{R}^3$ alternatives with vectors of attributes equidistant, in terms of ${\cal F}^{GL}(y)$ with $-\frac{1}{3}$, to $t\vec{1}, \ t \geq 0$. All such triangles form the surface of a wedge.}
	\label{Figure5A}
\end{figure}

A modification of (\ref{augLeontief})
\begin{equation}
	\label{eq4.2}
{\cal F}(y) = \min \{ \min_{l, \, 1 \leq l \leq k} [y_l + \rho \sum^k_{l=1}y_l] \, ;b\sum^k_{l=1}y_l \}, 
\end{equation}
for some $b > 0$ (Figure \ref{Figure6mat}), is also a PDIA function that is a locally a PDCA function.

Observe that a similar but smooth function (Figure \ref{Figure8mat}):
\begin{equation}
	\label{eq4.3}
{\cal F}(y) = -((\sum^k_{l=1}(\hat{y}_l - y_l)^p)^{1/p}  +\rho\sum^k_{l=1}(\hat{y}_l - y_l)), 
\end{equation}
where $1 \leq p \leq +\infty$, $\hat{y}_l = \max_{y \in Z}y_l(x), \ l = 1,\dots,k$, is a PDIA function but locally it is a PDCA function. 

All functions (\ref{Leontief_function}), 
 (\ref{Chebyshev}), (\ref{augLeontief}), (\ref{eq4.2}), (\ref{eq4.3}),
are symmetric to the half line $y = \vec{1}t, \ t \in [0:+\infty)$, because attribute transformation absorbing weights, as mentioned earlier. 

\begin{figure}
	\centering
		\subfigure[]{\includegraphics[width=0.45\textwidth]{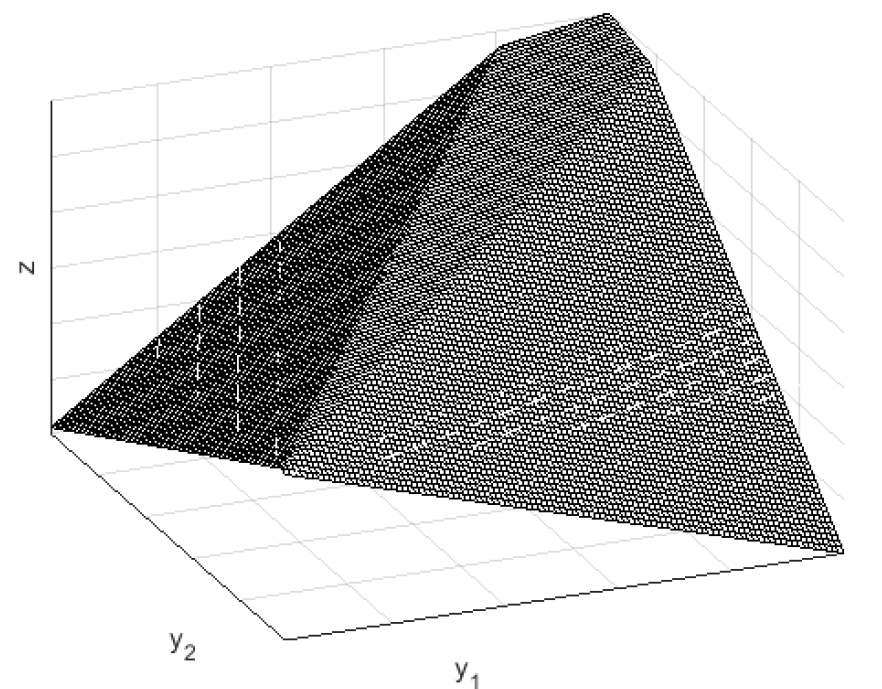}}
		\subfigure[]{\includegraphics[width=0.45\textwidth]{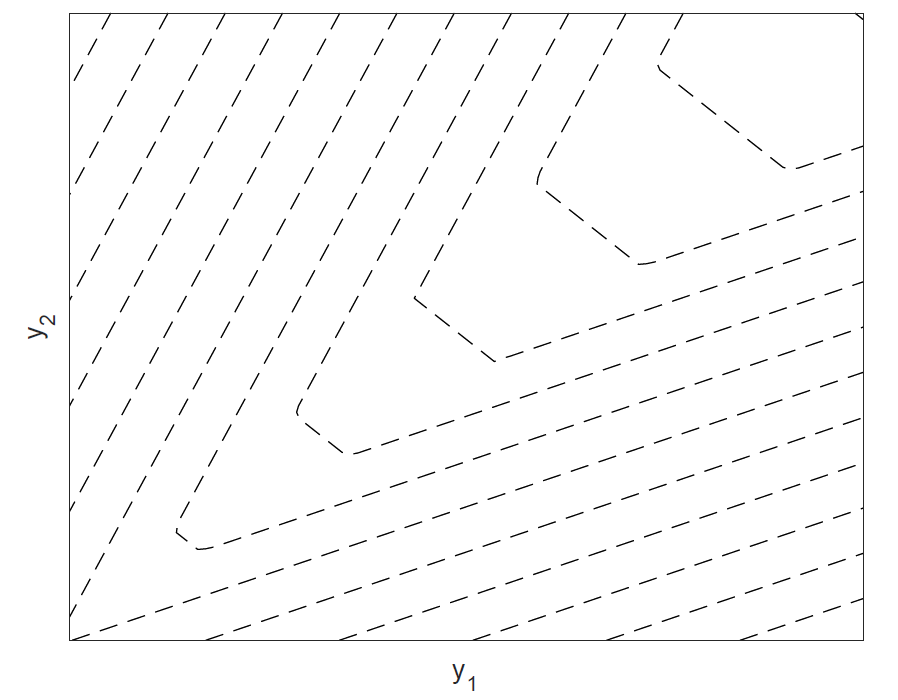}}
	\caption{Graph and contours of the PDIA function (\ref{eq4.2}) that is locally a PDCA function.}
	\label{Figure7mat}
\end{figure}

\begin{figure}
	\centering
		\subfigure[]{\includegraphics[width=0.45\textwidth]{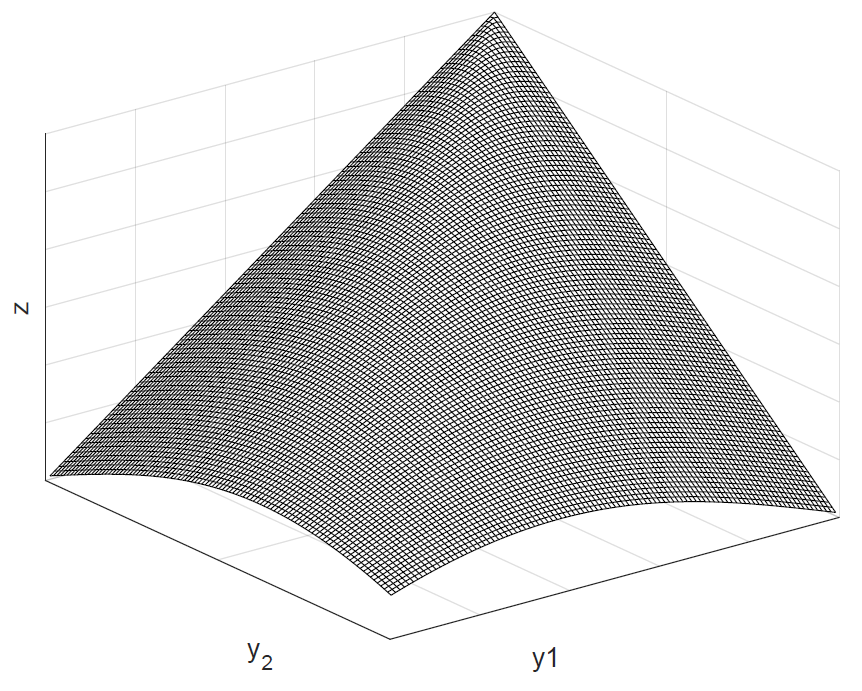}}
		\subfigure[]{\includegraphics[width=0.45\textwidth]{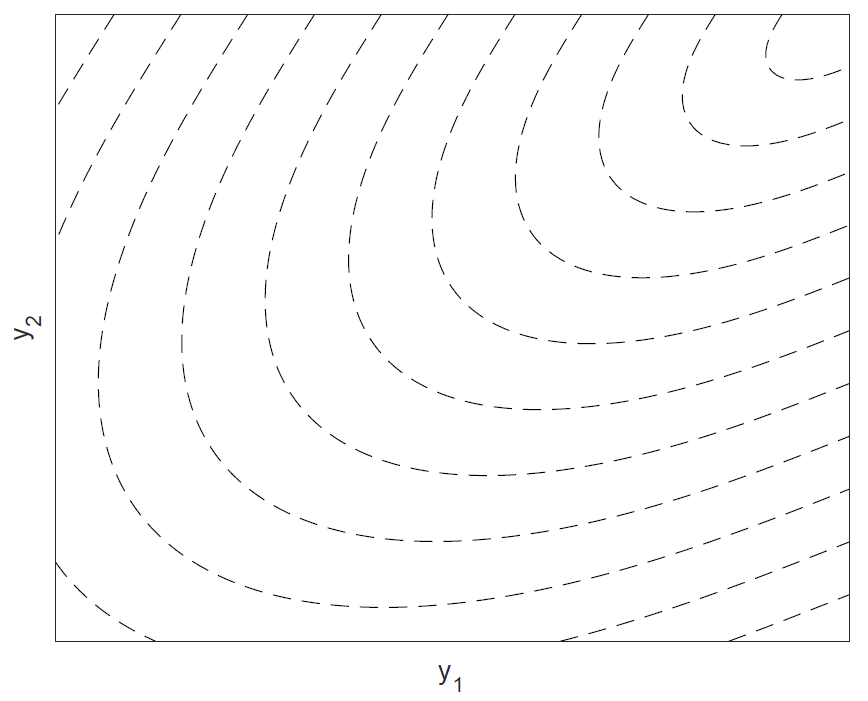}}
	\caption{Graph and contours of the smooth PDIA function (\ref{eq4.3}) that is locally a PDCA function.} 
	\label{Figure8mat}
\end{figure}

\section{Rankings with Pareto inconsistent assessment functions}
\label{nP_rankings}

We demonstrate how $K^{Poly}_{\rho}$ cones with $0 > \rho > -\frac{1}{k}$ and PDIA functions can be applied to a practical problem, namely to rankings.

Rankings are ubiquitous, and they are responsible for life-shaping decisions in the economics, technology, social and private lives. At present, rankings are produced with the help of PDCA functions. Thus they do not take much into account the attribute balance. However, in practical reality, alternatives with unbalanced attributes, assessed equally or similarly with alternatives with more balanced attributes can be undesired.
If so, PDIA functions should be applied instead.

\begin{figure}
	\centering
	\includegraphics[width=0.5\textwidth]{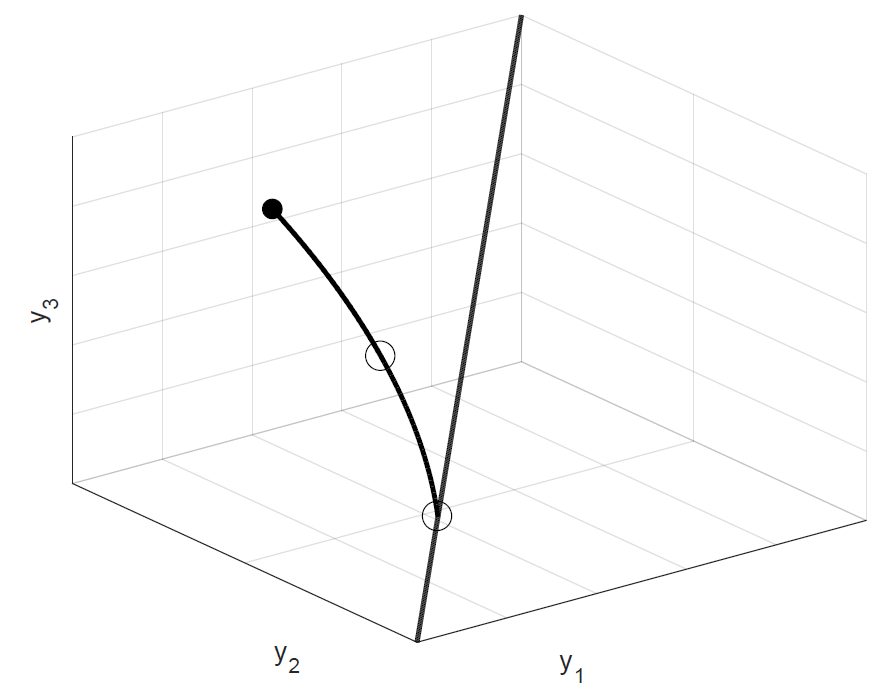}
	\caption{Example 1. All curve elements, except that marked by a bullet, are $R^k_+$-dominated.} 
	\label{Figure9}
\end{figure}

To illustrate the soundness of this idea,
we use data from the Environmental Performance Index (Environmental Performance Index 2023). EPI is calculated yearly for each country in the world (180 countries), collected, and made annually available by Yale and Colombia Universities (Wolf et al. 2022). 

EPI is the weighted sum over forty environmental performance indicators (i.e., country attributes), cumulated into three categories: climate change mitigation (acronym: PCC, cumulative weight 0.38), environmental health (acronym: HLT, cumulative weight 0.20), and ecosystem vitality (acronym: ECO, cumulative weight 0.42). So EPI is calculated as
\[
EPI = 0.38 PCC + 0.20 HLT + 0.42 ECO ,
\]
hence the suggested proportions between indicators PCC, HLT, and ECO is $(0.38:0.20:0.42)$.   

With $PCC' = 0.38 PCC$, $ HLT' = 0.20 HLT$, $ECO' = 0.42 ECO$, we have 
\[
EPI = PCC' + HLT' + ECO' ,
\]
and in consequence, we get $1:1:1$ as the suggested proportions between indicators PCC', HLT', and ECO' (a \textit{golden standard}). 

\begin{figure}
	\centering
	\includegraphics[width=0.5\textwidth]{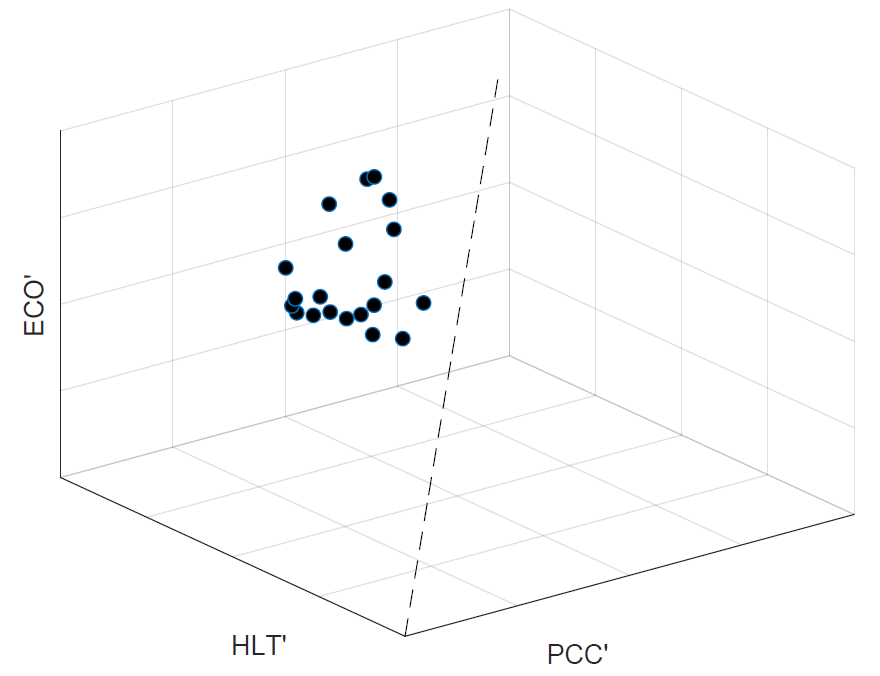}
	\caption{Example 2. The top-twenty EPI 2022 ranked countries represented by their performance indicators.} 
	\label{EPI_20_countries_3D}
\end{figure}

Figure \ref{EPI_20_countries_3D} presents the values of PCC', HLT', and ECO' and the values of respective assessment functions for the first twenty EPI 2022 ranked countries. It is interesting to note that the majority of countries are located well off the perfect balance half line $\vec{1}t, \ t \in [0:+\infty)$.

Table \ref{table1} presents the top-twenty EPI 2022 ranked countries and the top-twenty countries in three rankings according to (\ref{augLeontief}):  ranking for $\rho = 0$ (the Leontief ranking), the ranking for $\rho = -0.25 $, and, what we call the Perfect Balance Ranking (PBR), the ranking for $\rho = -\frac{1}{3}$ (only the top twenty countries in each ranking are presented). As $\rho$ decreases, the successive rankings more and more promote countries with more balanced attributes to higher positions and disfavor countries with lesser attribute balance.

\begin{table}
	\caption{The top-twenty countries according to four rankings}
	\label{table1}
	\centering   
	\begin{tabular}{cllll}
		& \multicolumn{4}{c}{Ranking} \\ \cmidrule{2-5}
		& \multicolumn{1}{c}{EPI} & \multicolumn{1}{c}{Leontief} & \multicolumn{1}{c}{ --} & \multicolumn{1}{c}{PBR}  \\ 
		&     &\multicolumn{1}{c}{$\rho= 0$} &\multicolumn{1}{c}{$\rho =-0.25$} & \multicolumn{1}{c}{$\rho=-\frac{1}{3}$} \\  \cmidrule{2-5}
		1	&	Denmark	&	Iceland	&	Ireland	&	T\"{u}rkiye	\\
		2	&	United Kingdom	&	Finland	&	Iceland	&	Israel	\\
		3	&	Finland	&	Sweden	&	Israel	&	Ireland	\\
		4	&	Malta	&	Switzerland	&	Singapore	&	Singapore	\\
		5	&	Sweden	&	Ireland	&	Brunei &	Brunei 	\\
		6	&	Luxembourg	&	Luxembourg	&	Norway	&	Uruguay	\\
		7	&	Slovenia	&	Denmark	&	Portugal	&	Algeria	\\
		8	&	Austria	&	United Kingdom	&	Australia	&	Iceland	\\
		9	&	Switzerland	&	France	&	T\"{u}rkiye	&	Malaysia	\\
		10	&	Iceland	&	Norway	&	Uruguay	&	Myanmar	\\
		11	&	Netherlands	&	Netherlands	&	Spain	&	Lebanon	\\
		12	&	France	&	Australia	&	USA	&	Papua \\
		13	&	Germany	&	Germany	&	Japan	&	Russia	\\
		14	&	Estonia	&	Austria	&	Switzerland	&	Oman	\\
		15	&	Latvia	&	Japan	&	New Zealand	&	Argentina	\\
		16	&	Croatia	&	Spain	&	France	&	Portugal	\\
		17	&	Australia	&	Belgium	&	Belgium	&	Maldives	\\
		18	&	Slovakia	&	Singapore	&	Netherlands	&	Qatar	\\
		19	&	Czech Republic	&	Italy	&	Argentina	&	Viet Nam	\\
		20	&	Norway	&	New Zealand	&	Italy	&	USA \\ 
	\end{tabular}
\end{table}

\section{Discussion and concluding remarks}
\label{conclusions}

The idea of PDIA functions has steamed independently from two domains that use quantitative multiattribute evaluations: the theory of production (Gadomski 2007) and the theory of multiple criteria decision-making (Kaliszewski 2006,2024). It undermines the well-established concept of Pareto efficiency that is considered a cornerstone of rational human actors, individual or collective. However, that concept is, as shown in this work, too restrictive to take into account the aspect of the attribute balance. 

The framework of PDIA functions developed here  is general enough to be applied, in principle, to rankings in any domain of human activities, whenever the attribute balance of perspective alternatives is an issue.

\setcounter{section}{1}
\renewcommand{\thesection}{\Alph{section}.2}
\section*{Appendix 1}
\label{appendix2}

We show that for a sequence of elements on the circle $y^2_1 + y^2_2 = 1$, tending to element $\bar{y} = (\frac{2\sqrt{5}}{5}, \ \frac{\sqrt{5}}{5})$ from the right side (i.e., with decreasing values of $y_1$), the sequence of corresponding numbers $\frac{r_2 - \bar{r}_2}{\bar{r}_1 - r_1}$ is unbounded.

Indeed,  
\[
\bar{r}_1 = 2\bar{y}_1 + \bar{y}_2 = 2\frac{2\sqrt{5}}{5} + \frac{\sqrt{5}}{5} = \sqrt{5}.
\] 
\[
\bar{r}_2 = \bar{y}_1 + 2\bar{y}_2 = \frac{2\sqrt{5}}{5} + 2\frac{\sqrt{5}}{5}= \frac{4\sqrt{5}}{5}.
\]
\[
\frac{y_1 + 2y_2 - \frac{4\sqrt{5}}{5}}{\sqrt{5} - 2y_1 - y_2}  = \frac{\sqrt{1 - y^2_2} + 2y_2 - \frac{4\sqrt{5}}{5}}{\sqrt{5} - 2\sqrt{1 - y^2_2} - y_2}
\]
At $\bar{y}_2 = \frac{\sqrt{5}}{5}$, the expression 
\[
lim_{y_2 \rightarrow^+ \frac{\sqrt{5}}{5}} \frac{\sqrt{1 - y^2_2} + 2y_2 - \frac{4\sqrt{5}}{5}}{\sqrt{5} - 2\sqrt{1 - y^2_2} - y_2} \,.
\]
is of the type $\frac{0}{0}$, so we apply L'H\^opital's rule, namely

\[
lim_{y_2 \rightarrow^+ \frac{\sqrt{5}}{5}} \frac{(\sqrt{1 - y^2_2} + 2y_2 - \frac{4\sqrt{5}}{5})'}{(\sqrt{5} - 2\sqrt{1 - y^2_2} - y_2)'} = lim_{y_2 \rightarrow^+ \frac{\sqrt{5}}{5}} \frac{-y_2 + 2\sqrt{1 -y^2_2}}{2y_2 - \sqrt{1-y^2_2}}.  
\]
When $y_2 \rightarrow^+ \frac{\sqrt{5}}{5}$ the numerator tends to the constant value $\frac{3\sqrt{5}}{5}$ whereas the denominator tends to zero. Thus, the sequence $\frac{r_2 - \bar{r}_2}{\bar{r}_1 - r_1}$ is unbounded.

\setcounter{section}{1}
\renewcommand{\thesection}{\Alph{section}.1}
\section*{Appendix 2}
\label{appendix4}

\noindent
{\it Proof of Theorem \ref{theorem4.1}}

The proof of this theorem is a generalization of the proof of Theorem \ref{theorem3.1}, given in Kaliszewski 1994, cf. also Steuer, Choo 1983, Choo, Atkins 1983. 

\noindent
{\bf Proof.}

\noindent
({\em Necessity})
Let $\bar{y} \in Z^\bullet$ be properly $K^{Poly}_{\rho}$-efficient. Then, there exists $N > 0$ such that for each
$\ l, \ l=1,\dots,k$, we have $r_{l} - \bar{r}_{l} \leq N(\bar{r}_{j}-r_{j})$
for some $j$
such   that   $r_{j}<\bar{r}_{j} \, ,$   whenever
$y\in Z^\bullet$   and  $r_{l}>\bar{r}_{l} \, .$
Let  $0 < \sigma \leq(N(k-1))^{-1} \,.$

Let $y$ be any element of $Z^\bullet \, , \ y \neq \bar{y} \, .$

\noindent
A. Suppose first that $\sum^k_{i=1} \bar{r}_i \geq \sum^k_{i=1} r_i $. 
Since $\bar{y}$ is $K^{Poly}_{\rho}$-efficient, there exists an index $h$
such that $\bar{r}_{h} > r_{h}$. Then,
\[
\bar{r}_h + \sum^k_{i=1} \bar{r}_i > r_h + \sum^k_{i=1} r_i.
\]
Since  $\bar{y}, y \in Z^{\bullet}$, 
\[
\bar{r}_h + \sum^k_{i=1} \bar{r}_i > r_h + \sum^k_{i=1} r_i > 0,
\]
and the last inequality implies that
\[
1 = \min_{l}\bar{\lambda}_{l}(\bar{r}_{l} + \sigma \sum^k_{i=1} \bar{r}_i) > 
\min_{l}\bar{\lambda}_{l}(r_{l} +
\sigma \sum^k_{i=1} r_i) \, ,
\]
where
\[\bar{\lambda}_{l} = (\bar{r}_{l} +
\sigma \sum^k_{i=1} \bar{r}_i)^{-1} > 0 \, , \ l=1,\dots,k . 
\]

\noindent
B. Suppose now
$\sum^k_{i=1}\bar{r}_i < \sum^k_{i=1} r_i $. 
Let now index $h$ be such that
$\bar{r}_{h} - r_{h} = \max_{i}(\bar{r}_{i} - r_{i}) \, .$ Since $\bar{y}$ is $K^{Poly}_\rho$-efficient, $\bar{r}_h - r_h > 0$. Moreover, since
$\bar{y}$ is properly $K^{Poly}_{\rho}$-efficient,  for
all $i$ such that
$r_{i} - \bar{r}_{i} > 0$, we have
\[(r_{i} - \bar{r}_{i})(\bar{r}_{h} - r_{h})^{-1} \leq  N,
\]
which entails
\[
0 < (\bar{r}_{h} - r_{h})^{-1} \sum_{i;r_{i} - \bar{r}_{i} > 0}
(r_{i} - \bar{r}_{i}) \leq N(k-1) \, ,
\]
and
\[
(N(k-1))^{-1}
\leq (\bar{r}_{h} - r_{h}) (\sum_{i;r_{i} - \bar{r}_{i} > 0}
(r_{i} - \bar{r}_{i}))^{-1} \,.
\] 
By the assumption, 
\[
0 < \sigma \leq (N(k-1))^{-1}
	\leq (\bar{r}_{h} - r_{h}) (\sum_{i;r_{i} - \bar{r}_{i} > 0}
	(r_{i} - \bar{r}_{i}))^{-1} .
\]

Since $\sum^k_{i=1}\bar{r}_i < \sum^k_{i=1} r_i $, 
we have $\sum^k_{i=1} (r_i - \bar{r}_i) > 0$ and 
\[
0 < \sigma \leq (N(k-1))^{-1}
\leq (\bar{r}_{h} - r_{h}) (\sum_{i;r_{i} - \bar{r}_{i} > 0}
(r_{i} - \bar{r}_{i}))^{-1} 
\]
\[
< (\bar{r}_{h} - r_{h})(\sum^k_{i=1} (r_i - \bar{r}_i))^{-1} \,  .
\]
Then,
\[\sigma \sum^k_{i=1} (r_i - \bar{r}_i) < (\bar{r}_{h} - r_{h})
\]
and
\[
\bar{r}_{h} + 
\sigma \sum^k_{i=1} \bar{r}_i > 
r_{h} + \sigma \sum^k_{i=1} r_i \, .
\]
Since  $\bar{y}, y \in Z^{\bullet}$,
\[
\bar{r}_{h} + 
\sigma \sum^k_{i=1} \bar{r}_i > 
r_{h} + \sigma \sum^k_{i=1} r_i > 0 \, .
\]
By the definition of
$\bar{\lambda} \, ,$ the last inequality implies that
\[
1 = \min_{l}\bar{\lambda}_{l}(\bar{r}_{l} + 
\sigma \sum^k_{l=1} \bar{r}) > \min_{l}\bar{\lambda}_{l}
( r_{l} + \sigma \sum^k_{l=1} r_l) \, .
\]

Since A or B holds for arbitrary $y \in Z, \ y \neq \bar{y}$, $\bar{y}$ solves $P^{+\infty}_{K^{Poly}_{\rho}}$ uniquely.

\vspace{0.2cm}
({\em Sufficiency})   Let
$\bar{y}  \in Z^\bullet$ solve the problem $P^{+\infty}_{K^{Poly}_{\rho}}$ for some
$\lambda  \in R^k_>$  and some $\sigma > 0$.  Suppose
$\bar{y}$ is not $K^{Poly}_{\rho}$-efficient. Then, for some
$y\in Z^\bullet \, , \ y \neq \bar{y}\, ,$ we have $r_l \geq \bar{r}_l, \ l=1,\dots,k$, and 
$r_h > \bar{r}_h $ for some index $h$.

Thus
\[
\sigma \sum^k_{i=1} r_i > \sigma \sum^k_{i=1}  \bar{r}_i \, ,
\]
and
\[
(r_{l} + \sigma \sum^k_{i=1} r_i) >
(\bar{r}_{l} + \sigma \sum^k_{i=1} \bar{r}_i) \, ,
\ l=1,\dots,k \, .
\]
Since  $\bar{y}, y \in Z^{\bullet}$,
\[
(r_{l} + \sigma \sum^k_{i=1} r_i) >
(\bar{r}_{l} + \sigma \sum^k_{i=1} \bar{r}_i) > 0 \, ,
\ l=1,\dots,k \, ,
\]
Since $\lambda \in R^k_>$,
\[
\lambda_{l}(r_{l} + \sigma \sum^k_{i=1} r_i) >
\lambda_{l}(\bar{r}_{l} + \sigma \sum^k_{i=1} \bar{r}_i) \, ,
\ l=1,\dots,k \, ,
\]
and finally
\[
\min_{l}\lambda_{l}(r_{l} + \sigma \sum^k_{i=1} r_i) >
\min_{l}\lambda_{l} (\bar{r}_{l} +
\sigma \sum^k_{i=1} \bar{r}_i) \, ,
\]
but this is a contradiction because $\bar{y}$ solves problem $P^{+\infty}_{K^{Poly}_\rho}$. Hence,
$\bar{y}$ is $K^{Poly}_{\rho}$-efficient.

Denote $\tau_{l}= r_{l} + \sigma \sum^k_{i=1} r_i$  and $\bar{\tau}_{l} = \bar{r}_{l} + \sigma \sum^k_{i=1}  \bar{r}_i$, $l =1,\dots,k$. 
Observe that
\[
\tau_{l} - \bar{\tau}_{l} = \sum_{i\neq l}\sigma (r_{i} - \bar{r}_{i}) + (\sigma + 1)
(r_{l} - \bar{r}_{l}) , 
\]
for all $l = 1,...,k$.

By $K^{Poly}_{\rho}$-efficiency of $\bar{y}\, ,$ for any $y\in Z^\bullet , \ y\neq \bar{y}\, ,$ there exists an index $h$ such that  $r_{h} < \bar{r}_{h}$\,.
Let $\tau_l - \bar{\tau}_l = \min_i( \tau_i - \bar{\tau}_i)$. From  the
efficiency of  $\bar{r}$  it follows that
$\tau_{l} - \bar{\tau}_{l} \leq 0$ and this implies
$r_{l}-\bar{r}_{l}<0 \, .$ In fact, suppose
$r_{l}- \bar{r}_{l}\geq 0 \, .$  Then
\[
\tau_{h} - \bar{\tau}_{h}\begin{array}[t]{ll}
	=(r_{h} - \bar{r}_{h}) + \rho e^{k}( r - \bar{r}) \\
	<(r_{l} - \bar{r}_{l}) + \rho e^{k}( r - \bar{r}) \\
	= \tau_{l} - \bar{\tau}_{l} \\
	=\min_{i}(\tau_{i} - \bar{\tau}_{i}) \, ,
\end{array}
\]
but this is a contradiction. 

Moreover,
$\min_{i}(\tau_{i} - \bar{\tau}_{i}) =\tau_{l} - \bar{\tau}_{l} $ implies that 
\[
\tau_{l} - \bar{\tau}_{l} \leq \tau_{i} - \bar{\tau}_{i}, \ i = 1,...,k,
\]
\[
(r_{l} - \bar{r}_{l}) + \rho e^{k}(\bar{r} - r) \leq (r_{i} - \bar{r}_{i}) + \rho e^{k}(\bar{r} - r), \ l = 1,...,k,
\]
\[
r_{l} - \bar{r}_{l} \leq r_{i} - \bar{r}_{i}, \ l = 1,...,k,
\]
\[r_{l} - \bar{r}_{l} = \min_{i}(r_{i} - \bar{r}_{i} ).
\]

Let
$r_{j} - \bar{r}_{j} = \max_{i}(r_{i} - \bar{r}_{i})$
and  suppose
$r_{j} - \bar{r}_{j}>0 \,.$  
Then, 
\[
0 > \tau_{l} - \bar{\tau}_{l} = \sum_{i\neq l}\sigma (r_{i} - \bar{r}_{i}) + (\sigma + 1)
(r_{l} - \bar{r}_{l})
\]
\[  
	\geq (k-2)\sigma (r_{l} - \bar{r}_{l}) + \sigma
	(r_{j} - \bar{r}_{j}) + (\sigma +1)( r_{l} - \bar{r}_{l}) 
\]
\[
	= (1 + (k-1)\sigma)( r_{l} - \bar{r}_{l}) + \sigma(r_{j} - \bar{r}_{j})
\]
and
\[
\frac{r_{j} - \bar{r}_{j}}{\bar{r}_{l} - r_{l}} \leq \frac{1+(k-1)\sigma}{\sigma} \, .
\]

Finally, $N \leq (1+(k-1)\sigma))\sigma^{-1} \, .$ \mybox

\setcounter{section}{1}
\renewcommand{\thesection}{\Alph{section}.3}
\section*{Appendix 3}
\label{appendix3}

\setcounter{Lemma}{0}
\begin{Lemma}
	For $\rho = -\frac{1}{k}$,
		contours of function (\ref{GLF}) reduce to half lines parallel to the half line $t \vec{1}, \ t \in [0,+\infty)$. For any point on the half line $t \vec{1}, \ t \in [0,+\infty)$, the function assumes value $0$.
\end{Lemma}

\vspace{0.2cm} 
\noindent 
{\bf Proof.} 
For any half line $\bar{y} + t\vec{1}, \ t \in [0,+\infty), \ \bar{y} \in \mathbb{R}^k$, parallel to the half line $t\vec{1}, \ t \in [0,+\infty)$, consider an element $y$ on that line, $y = \bar{y} + t\vec{1}, t \in [0:+\infty)$.
Then we have
\[
{\cal F}(y) = \min_{l; \, 1 \leq l \leq k} (y_l - \frac{1}{k} \sum^k_{i=1} y_i) = \min_{l; \, 1 \leq l \leq k} (\bar{y}_l + t - \frac{1}{k} \sum^k_{i=1 }(\bar{y}_i + t) ) 
\]
\[
= \min_{l; \, 1 \leq l \leq k} (\bar{y}_l - \frac{1}{k} \sum^k_{i=1 }\bar{y}_i + t - \frac{1}{k} \sum^k_{i=1 }t )
\]
\[
= \min_{l; \, 1 \leq l \leq k} (\bar{y}_l - \frac{1}{k} \sum^k_{i=1 }\bar{y}_i) , 
\]
that is a constant value irrespective of $t$. 

Any point $y$ on the half line $t \vec{1}, \ t \in [0:+\infty)$, has all its components equal, hence $\frac{1}{k} \sum^k_{l=1 } y_l = y_l$ for any $l$. Thus, the function (\ref{GLF}) for $\rho = -\frac{1}{k}$ assumes value $0$. \mybox

\end{document}